\documentclass[12pt,twoside]{article}
\usepackage{amsmath}
\usepackage{amsfonts}
\usepackage{amsthm}
\usepackage{amssymb}
\usepackage{graphicx}
\usepackage{multicol}

\begin{document}
\title{{\normalsize
{\bf Note on surface-link of trivial components}}}
\author{{\footnotesize Akio Kawauchi}\\
\date{}
{\footnotesize{\it Osaka Central Advanced Mathematical Institute, Osaka Metropolitan University}}\\ 
{\footnotesize{\it Sugimoto, Sumiyoshi-ku, Osaka 558-8585, Japan}}\\ 
{\footnotesize{\it kawauchi@omu..ac.jp}}}
\date\, 
\maketitle
\vspace{0.1in}
\baselineskip=9pt
\newtheorem{Theorem}{Theorem}[section]
\newtheorem{Conjecture}[Theorem]{Conjecture}
\newtheorem{Lemma}[Theorem]{Lemma}
\newtheorem{Sublemma}[Theorem]{Sublemma}
\newtheorem{Proposition}[Theorem]{Proposition}
\newtheorem{Corollary}[Theorem]{Corollary}
\newtheorem{Claim}[Theorem]{Claim}
\newtheorem{Definition}[Theorem]{Definition}
\newtheorem{Example}[Theorem]{Example}

\begin{abstract} As a previous result, it has shown that every sphere-link consisting of trivial components is a ribbon sphere-link. 
In this note, it is shown that for every closed oriented disconnected surface ${\mathbf F}$ 
with just one non-sphere component, every ${\mathbf F}$-link consisting of trivial components 
is a ribbon surface-link. 
Further, it is shown that for every closed oriented disconnected surface ${\mathbf F}$ 
containing at least two non-sphere components, there exist a pair of a ribbon 
${\mathbf F}$-link and 
a non-ribbon ${\mathbf F}$-link that consist of trivial components and have meridian-preservingly 
isomorphic fundamental groups.

\medskip

\noindent{\it Keywords}: Ribbon,\, Non-ribbon, Surface-link.

\noindent{\it Mathematics Subject Classification 2010}: Primary 57Q45; Secondary 57N13
\end{abstract}

\baselineskip=15pt

Let ${\mathbf F}$ be a  (possibly disconnected) closed surface. 
An ${\mathbf F}$-{\it link}  in the 4-sphere $S^4$ is the image of  a smooth embedding 
${\mathbf F}\to S^4$. When ${\mathbf F}$ is connected, it is also called 
an ${\mathbf F}$-{\it knot}. An ${\mathbf F}$-link or 
${\mathbf F}$-knot for an ${\mathbf F}$ is called a  {\it surface-link} or {\it surface-knot} in $S^4$, respectively. 
If ${\mathbf F}$ consists of some copies of the 2-sphere $S^2$, then it is also called an 
$S^2$-{\it link}. and an $S^2$-{\it knot} for ${\mathbf F}=S^2$.
A {\it trivial  surface-link} is a surface-link $F$ which bounds disjoint 
handlebodies smoothly embedded in $S^4$. 
A {\it ribbon surface-link} is a surface-link $F$ which is obtained 
from a trivial $n S^2$-link $O$ for some $n$ (where $n S^2$ denotes 
the disjoint union of $n$ copies of  $S^2$) 
by surgery along an embedded 1-handle system, \cite{K15}, \cite{K17-1}, \cite{K17-2}, \cite{K18},  \cite{KSS}. 
It was shown that every $S^2$-link consisting of trivial components is a ribbon 
$S^2$-link, \cite{K24}. 
The following result is a generalization of this result. 

\medskip

\medskip

\noindent{\bf Theorem~1.} 
Let ${\mathbf F}$ be a closed oriented disconnected surface with at most one non-sphere component. Then 
every ${\mathbf F}$-link $L$ in $S^4$ consisting of trivial components is a ribbon 
${\mathbf F}$-link in $S^4$.

\medskip 

\medskip

\noindent{\it Proof of Theorem~1.} The case that ${\mathbf F}$ consists of only $S^2$-components has been given, \cite[Theorem~1.5]{K24}. This proof is done by a similar method. 
Let ${\mathbf F}$ have $S^2$-components and 
only one non-sphere component, and $L$ an ${\mathbf F}$-link in $S^4$ consisting of trivial components. 
Let $F$ be the trivial non-sphere component of $L$   and $L'=L\setminus F$ the remaining 
$S^2$-link consists of trivial components. Since the second homology class $[F]=0$ in 
$H_2(S^4\setminus L';Z)=0$, there is a compact connected oriented 3-manifold 
$V_F$ smoothly embedded in $S^4$ with $\partial V_F=F$ and $V_F\cap L'=\emptyset$. 
The $S^2$-link $L'$ is a ribbon $S^2$-link in $S^4$, \cite[Theorem~1.5]{K24}. 
Let $W$ be a SUPH system in $ S^4$ with $\partial W=L'\cup O$
for a trivial $S^2$-link $O$, \cite{K24}. 
Let $\alpha$ be an arc system in $W$ spanning $O$ such that the closed complement 
$\mbox{cl}(W\setminus N(\alpha))$ is identified with a boundary collar $L'\times[0,1]$ 
of $L'$ in $W$ with $L'\times 0=L'$ where $N(\alpha)$ is a regular neighborhood of $\alpha$ in $W$ which is written as a trivial disk fiber bundle $d\times \alpha$ over $\alpha$. 
Since $V_F\cap L'=\emptyset$, it can be assumed that $(L'\times[0,1])\cap V_F=\emptyset$ 
and the interior of the disk fiber $d\times x$ for a point $x$ of $\alpha$ meets $V_F$ transversely with a simple proper arc system and a simple loop system. By deforming $V_F$, the interior of the disk fiber $d\times x$ meets $V_F$ transversely only with a simple proper arc system $\beta$. As a result, the intersection $N(\alpha)\cap V_F$ is assumed to be a thickening $\beta\times [0,1]$ of 
$\beta$. A regular neighborhood of $\beta\times [0,1]$ in $V_F$ is a 
1-handl system $h(\beta)$ on $F$. By adding a disjoint 1-handle system $h^+$ on $F$ embedded in $V_F$, the closed complement $H=\mbox{cl}(V_F\setminus (h(\beta)\cup h^+))$
is a handlebody of genus, say $n$. Let $F^+=\partial H$.
Let $H^0$ be a once-punctured handlebody obtained from $H$ by removing a 3-ball 
with $\partial H^0=F^+\cup O^H$. 
The union $W\cup H^0$ is a SUPH system for the surface-link $L'\cup F^+$ in $S^4$, so that $L'\cup F^+$ is a ribbon surface-link in $S^4$ by \cite[Lemma~3.1]{K24}. 
The 1-handle system $h(\beta)\cup h^+$ on $F$ is a 2-handle system 
$D_i\times  I \, (i=1,2, \dots   , n)$ on $F^+$, where $I$ denotes an interval containing $0$ in the interior. Then there is a disjoint O2-handle pair system 
$(D_i \times I, D'_i\times I) \, (i = 1, 2, \dots, n)$ on $F^+$ because any disjoint 1-handle system on the trivial surface-knot $F$ is a disjoint trivial 1-handle system on $F$, \cite{HoK},  
\cite[Lemma 4.1]{K24}. Let $(O\cup O^H, \alpha\cup \alpha^H)$  be a chorded sphere system for the ribbon surface-link  $L'\cup F^+$ constructed in the SUPH system $W\cup H^0$. Let $B(O)\cup B(O^H)$ be a disjoint 3-ball system bounded by the trivial $S^2$-link $O\cup O^H$  
in $S^4$. The intersections $B(O)\cap D_i = B(O)\cap D'_i = \emptyset (i =1, 2, \dots   , n)$  can be assumed by moving the 3-ball system $B(O)\cup B(O^H)$ in $S^4$. The intersections 
$(\alpha\cup \alpha^H)\cap D_i 
= (\alpha\cup \alpha^H)\cap D'_i
= \emptyset\, (i =1, 2, \dots , n) $ 
can be also assumed by general position. Thus, the disjoint O2-handle pair system 
$(D_i \times I, D'_i\times  I) \, (i = 1, 2, \dots   , n)$ on $F$ is deformed into a disjoint O2-handle pair system on the ribbon surface-link $L'\cup F^+$ in $S^4$, whose surgery surface-link $L = L'\cup F$ is a ribbon surface-link, 
\cite[ Corollary 1.2, \& Lemma 2.1]{K24}. 
This completes the proof of Theorem~1. 

\medskip

\medskip

In the case that ${\mathbf F}$ has at least two non-sphere components, the following result is obtained.

\medskip

\medskip

\noindent{\bf Theorem~2.} 
Let ${\mathbf F}$ be any closed oriented disconnected surface with at least two non-sphere components. Then there exist a pair of a ribbon ${\mathbf F}$-link $L$ and a non-ribbon 
${\mathbf F}$-link $L'$ in $S^4$ that consist of trivial components and have meridian-preservingly 
isomorphic fundamental groups.

\medskip

\medskip

\noindent{\it Proof of Theorem~2.} 
 Let $k\cup k'$ be a non-splitable link in the interior of  a 3-ball $B$ such that $k$ and $k'$ are trivial knots. 
For the boundary 2-sphere $S=\partial B$ and the disk $D^2$ with the boundary circle $S^1$, 
let $L$ be the torus-link consisting of the torus-components $T= k\times S^1$ 
and $T'= k'\times S^1$ in the 4-sphere $S^4$ with 
$S^4=B\times S^1 \cup S\times D^2$, which is a ribbon torus-link in $S^4$, \cite{K15}. 
Since $k$ and $k'$ are trivial knots in $B$, the torus-knots $T$ and $T'$ are trivial torus-knots in $S^4$ by construction. 
Since $k\cup k'$ is non-splitable in $B$, there is a simple loop $t(k)$ in $T$ coming from the longitude of $k$ in $B$ such that $t(k)$ does not bound any disk not meeting $T'$ in $S^4$, 
meaning that there is a simple loop $c$ in $T$ unique up to isotopies of $T$ 
which bound a disk $d$ in $S^4$ not meeting $T'$, where $c$ and $d$ are given by $c=p\times S^1$ and $d=a\times S^1\cup q\times D^2$ for a simple arc $a$ in $B$ joining a point $p$ of $k$ to a point $q$ in $S$ with $a\cap(k\cup k')=\{p\}$ and $a\cap S=\{q\}$. 
Regard the 3-ball $B$ as the product $B=B_1\times [0,1]$ for a disk $B_1$. 
Let $\tau_1$ is a diffeomorphism of the solid torus $B_1\times S^1$ given by one full-twist along the meridian disk $B_1$, and $\tau=\tau_1\times 1$ the product diffeomorphism of 
$(B_1\times S^1)\times[0,1]=B\times S^1$. 
Let $\partial\tau$ be the diffeomorphism of the boundary $S\times S^1$ of $B\times S^1$ 
obtained from $\tau$ by restricting to the boundary, 
and the 4-manifold $M$ obtained from $B\times S^1$ and $S\times D^2$ by pasting the boundaries 
$\partial (B\times S^1)=S\times S^1$ and $\partial (S\times D^2)=S\times S^1$ by the diffeomorphism $\partial\tau$. Since the diffeomorphism $\partial\tau$ of $S\times S^1$ extends 
to the diffeomorphism $\tau$ of $B\times S^1$, the 4-manifold $M$ is diffeomorphic to  $S^4$. 
Let $L_M=T_M\cup T'_M$ be the torus-link in the 4-sphere $M$ arising from 
$L=T\cup T'$ in $B\times S^1$. 
The fundamental groups $\pi_1(S^4\setminus L, x)$ and $\pi_1(M\setminus L_M, x)$ 
are meridian-preservingly isomorphic by van Kampen theorem. 
The loop $t(k)$ in $T_M$ does not bound any disk not meeting $T'_M$ in $M$, 
so that the loop $c$ in $T_M$ is a unique simple loop up to isotopies of $T_M$ which 
bounds a disk $d_M=a\times S^1\cup D^2_M$ in $M$ not meeting $T'_M$, where 
$D^2_M$ denotes a proper disk in $S\times D^2$ bounded by the loop 
$\partial\tau(q\times S^1)$. 
 An important observation is that the self-intersection number  
$\mbox{Int}(d_M, d_M)$  in $M$ with respect to the surface-framing on $L_M$ is $\pm1$. This means that the loop $c$ in $T_M$ is a non-spin loop and thus, 
the torus-link $L_M$ in $M$ is not any ribbon torus-link, \cite{HiK}, \cite{K02}. 
Let $(S^4,L')=(M,L_M)$. 
If ${\mathbf F}$ consists of two tori, then the pair $(L,L')$ forms a desired pair. 
If ${\mathbf F}$ is any surface consisting of two non-sphere components,   
then a desired ${\mathbf F}$-link pair is obtained from the pair $(L,L')$ by taking connected sums of some trivial surface-knots, because every stabilization of a ribbon
 surface-link is a ribbon surface-link and every stable-ribbon surface-link is a ribbon 
 surface-link, \cite{K24}. 
If ${\mathbf F}$ has some other surface ${\mathbf F}_1$  in addition to a surface
${\mathbf F}_0$ of two non-sphere components,
then  a desired ${\mathbf F}$-link pair is obtained from a desired ${\mathbf F}_0$-link pair  
by adding the trivial ${\mathbf F}_1$-link as a split sum. 
Thus, a desired ${\mathbf F}$-link pair $(L,L')$ is obtained. 
This completes the proof of Theorem~2. 

\medskip

\medskip

In the proof of Theorem~2, the diffeomorphism $\partial \tau$ of $S\times S^1$ coincides 
with  Gluck's non-spin diffeomorphism of $S^2\times S^1$, \cite{Gluck}. The torus-link $(M,T_M)$ called a {\it turned torus-link} of a link $k\cup k'$ in $B$ is an analogy of a turned torus-knot of a knot  in  $B$, \cite{Boyle}. There is an invariant of a surface-knot used to confirm a non-ribbon surface-knot, 
\cite{K02}. This invariant is easily generalized to an invariant of a surface-link, which  can be also applied to confirm that $L'$ is a non-ribbon surface-link.

\medskip

\medskip

\noindent{\bf Acknowledgements.} 
This work was partly supported by JSPS KAKENHI Grant Number JP21H00978 and MEXT Promotion of Distinctive Joint Research Center Program JPMXP0723833165.

\end{document}